*Draft of November 2021—please do not cite without permission.*

**A cautionary note about self-reference**[1]

T. Parent (Nazarbayev University)
nontology@gmail.com

## *1. Introduction*

Yablo (1993) argues that paradox can be generated without self-referential expressions. But he denies that "self-reference *suffices* for paradox" (p. 251). This is obvious if he means that not all cases of self-reference are paradoxical (e.g., 'This sentence is a sentence.'). Yet since he justifies the claim by citing Tarski and Gödel, Yablo likely meant something stronger. An alternate reading is: Self-reference does not suffice for *any* paradoxical sentence, at least in a semantically open language—that is to say, in a language without semantic terms ('true', 'refers, 'satisfies', etc.) which are defined on expressions of that self-same language.[2] This is a standard view; however, in what follows it is shown to be only approximately correct. The limitations of the view seem worth noting; nevertheless, we shall see that there are no revisionary consequences implied in this.

Caveat: Endorsements of unrestricted self-reference are hard to find in the literature. As far as I know, Gödel and Tarski never say as much; Yablo's remark is the closest I can find. More often, self-referential devices seem to enjoy a "presumption of innocence," as seen in the work of Tarski and Gödel. (Consider also Henkin's 1949 completeness proof for first-order logic.) At any rate, I expect that working logicians will recognize the presumption of innocence as orthodox; however, I shall suggest that some hedging is required.

---

[1] I thank Bradley Armour-Garb, John G. Bennett, Matteo Bianchetti, Ethan E. Brauer, Tim Button, David Chalmers, Morgan Davies, Anil Gupta, Thomas Hofweber, Cory Johnson, Yael Loewenstein, William Lycan, David McCarty, Jay Newhard, Bryan Pickel, Alexander Pruss, Michael Resnik, David Ripley, and Keith Simmons for discussion on various versions of this material. I also thank an audience at the 2018 Central Division meeting of the American Philosophical Association.

[2] The term 'semantically open language' is from Tarski (1944). Tarski restricted himself to such languages so to bypass Liar-like paradoxes, yet this is no longer popular due to alternate approaches by Kripke (1972) and Gupta & Belnap (1993) (although Kripke's is patently nonclassical given that it rejects the law of excluded middle). Still, it is widely assumed that *if* a language (/fragment) is semantically open, unrestricted self-reference alone is insufficient for paradox. This is what I shall examine in this paper.

*2. An argument for absurdity*

It is well known that substitution of co-referring terms is not truth preserving when the substitution occurs inside quotation marks. Suppose that the names or constants of the object language are '$a_0$', '$a_1$', '$a_2$',… And suppose that '$a_0$' and '$a_1$' are interpreted as co-referring, whence it is true that:

$$a_0 = a_1$$

Then, given the truism that '$a_1$' = '$a_1$', substitution inside quotation would allow us to deduce the absurdity that: '$a_0$' = '$a_1$'. Thus, the substitution into quotation marks is invalid.

What amounts to much the same, a *substitutable variable* inside quotation marks is sufficient to express an ill-defined function.[3] Thus, consider the expression:

'*x*'

This would normally mention the symbol used as a variable, viz., the (italicized) 24th letter of the alphabet. But suppose instead that the variable inside the quotes is substitutable, so that it is functionally equivalent to a "blank" that can be "filled in" by any term. Then, the quoted variable expresses an ill-defined function $q(x)$ = '$x$', one which yields different outputs for $a_0$ depending on which term for $a_0$ replaces the variable:

$q(a_0)$ = '$a_0$'

$q(a_1)$ = '$a_1$'

---

[3] This connection between substitution into quotes and an ill-defined function was observed already by, e.g., Tarski (1933/1983, p. 161).

In order to preclude this, we must insist that a quoted variable is always merely a name for the symbol enclosed in quotes, nothing more.[4]

Notwithstanding, an analogous ill-defined function can be re-created by exploiting the possibilities for self-reference, even in a semantically open language. Suppose here that we stipulate the notion of a *reflection* as follows. Where ‘$n$’ is a variable for a natural number with numeral $\underline{n}$:

(*) The reflection of $x = y$ iff there is an $n$ such that $x = a_n$ and $y = \ulcorner a_{\underline{n}} \urcorner$.

Notation: Per Quine (1951), an expression with corner quotes denotes the concatenation of symbols enclosed in the corners, after the replacement of any metavariables. Thus, as concerns $\ulcorner a_{\underline{n}} \urcorner$, if ‘$\underline{n}$’ is replaced with ‘0’, the corner-quoted expression then denotes the concatenation of the symbol ‘a’, followed by a subscripted ‘0’.[5] That is to say, it denotes ‘$a_0$’.

The definition at (*) institutes a type of self-reference—for when ‘$x$’ is replaced with a constant, (*) defines ‘the reflection of $x$’ to denote whichever constant is appended to that very descriptor. That is, when ‘$x$’ is replaced with a constant, the descriptor effectively means “The constant hereby concatenated with ‘the reflection of $x$’.”[6] If self-reference is unrestricted in a semantically open language, such a descriptor seems thus far admissible.

However, (*) introduces a function that is ill-defined on input $a_0$ in exactly the same way as $q(x)$. After all, (*) implies that:

[4] This should not suggest a problem with the function that maps an expression onto its quotation. If $\alpha$ is an expression, $f(\alpha) = \ulcorner‘\alpha’\urcorner$ remains perfectly well defined; consider that $f$(‘$a_0$’) = ‘‘$a_0$’’ and $f$(‘$a_1$’) = ‘‘$a_1$’’. (The outputs here are different, of course, but so are the inputs.)

[5] It may be desirable here and with similar sentences to use the notation of Boolos (1995), so to disambiguate which quotation marks are paired together.

[6] When ‘$x$’ is replaced with a saturated function symbol, however, ‘the reflection of $x$’ does not denote the very term appended to the descriptor. It rather denotes a constant that co-refers with the relevant functor. But such cases will not be relevant here.

(1) The reflection of $a_0$ = '$a_0$';

(2) The reflection of $a_1$ = '$a_1$'.

So like $q(x)$, the reflection "function" varies between outputting '$a_0$' and '$a_1$' on input $a_0$, depending on which name for the input used. Concurrently, whereas the definition at (*) rules that (1) is true, (1) is provably equivalent to the following absurdity:

(3) The reflection of $a_1$ = '$a_0$'.

The absurdity of (3) is even clearer in noting that, by the transitivity of identity, (2) and (3) imply that '$a_0$' = '$a_1$'. (N.B., the substitution needed to show the equivalence between (1) and (3) is *not* a substitution inside quotation marks.)

The upshot is that since ill-defined are to be excluded, the language must somehow be regulated to prevent the introduction of such a thing. So contra the orthodoxy, it is not strictly true that self-reference can be literally unrestricted in a semantically open language.

*3. Discussion*

Further elaboration is in order. First, note that the underline in (*) expressing the numeral function is not equivalent to quotation; witness that $\underline{0} = \underline{0+0}$ whereas '0' ≠ '0+0'. However, '$\underline{n}$' may appear to be a semantic expression, for $\underline{n}$ is the numeral that *denotes* the number *n*. If so, then the language may not be semantically open. However, the pathology ultimately depends on neither '$\underline{n}$', nor the underline, nor quotes, nor corner quotes. It is enough to have a (well-defined) function *h* such that:

$h(0)$ = the first lowercase letter of the alphabet with subscript nought;

$h(1)$ = the first lowercase letter of the alphabet with the successor numeral to nought as its subscript.

On all other inputs, assume *h* is undefined. We then can make do with a partial definition of a reflection:

(**) The reflection of $x = y$ if there is an $n$ such that $x = a_n$ and $y = h(n)$.

Assuming $a_0 = a_1$, the absurdity at (3) is derivable in a similar fashion.[7]

It is natural to think that the problem owes to (*) and (**) quantifying into the subscript position. But while this seems to be sufficient for the problem, I do not believe it is necessary (see the Appendix). Besides, classical logicians otherwise have need of quantification into subscript position. E.g., when introducing infinitely many constants into the language, strictly speaking one must invoke clauses such as:

- For all $n$, $\ulcorner a_{\underline{n}} \urcorner$ is a constant.
- For all $n$, $\ulcorner a_{\underline{n}} \urcorner$ denotes $a_n$.

Quantification into subscript position occurs elsewhere too, as in some presentations of the anti-diagonal function in Cantorian arguments.[8] (I do not mean to suggest a problem with such presentations; the point instead is that our argument requires only standard formal devices.)

It can be said that the reflection of *x* is an intensional notion, for the matter depends on what term for *x* is being used. This is noteworthy insofar as it is a distinctive kind of intensionality; it does not owe to a propositional attitude verb, an idiom like ‘so called’, attribute abstraction, a modal operator, or (to repeat) substitution inside quotes.[9] More important, however, is that the definition at (**) clearly uses only formal, nonsemantic terms. This reveals that we can prove an absurdity in a semantically open language simply by exploiting the possibilities for self-reference. Or from another angle, unchecked self-reference enables intensional contexts, even if the language lacks the usual intensional elements.

---

[7] My thanks to Alexander Pruss for raising the concern about ‘$\underline{n}$’ and suggesting a function like *h* in response.
[8] E.g., see the description of the anti-diagonal set of natural numbers in Boolos et al. (2007, ch. 2).
[9] On such cases of intensionality, see Quine (1960, sections 30ff.).

One may object that the preceding shows merely that definitions like (*) are not legitimate. But that is exactly the point: Such definitions must be excluded. Still, a definition like (*) seems problematic simply because 'The reflection of *x*' fails to express a well-defined function—and ill-defined functions are already excluded. Given that, one may insist that we never had need to worry about (*).

Basically, the point is correct; indeed, it explains why the argument does not have any drastic revisionary implications. At the same time, special attention is called for when commonplace formal devices give rise to an ill-defined function. This is how one might regard the ill-defined "quoting function" from section 2. Quotation seems unremarkable at first glance; it is just a convenient way to introduce metalinguistic terms. But substitution inside quotation makes for an ill-defined function, and so, a word of warning is required.

In the same way, self-reference has been "presumed innocent" by the working logician. But unrestricted self-reference is sufficient for an ill-defined function, even in a semantically open language—and so, an explicit advisory is apt. Put differently, the ban on ill-defined functions *implies* certain restrictions on self-reference; however, this seems not to be widely acknowledged.

*4. A version with Gödel numbers*

In conversation, Tim Button has worried that if unrestricted self-reference suffices for such a problem, then Peano Arithmetic might be shown unsound via the method of Gödel numbering. For Gödel numbering enables something functionally like self-reference, and unrestricted self-reference now appears sufficient for pathology.

A derivative argument is indeed readily available using Gödel numbers, as I shall explain shortly. However: The argument shows only that Gödel numbering in the metalanguage cannot be used with abandon. It hardly implies that the arithmetical object language is defective. (Indeed, nothing in Peano Arithmetic functions like self-reference prior to Gödel numbering.[10]) And again, the general prohibition on ill-defined functions means that, in one respect, there was never cause for alarm. Even so, it is worth remarking that the ban on ill-defined functions implies certain restraints on self-referential devices, including Gödel numbering. But to be absolutely clear: The Gödel numbering version of the "paradox" does not show that any use of Gödel numbering is bankrupt. (An analogy: The quoting function from section 2 does not show all use of quotation to be illegitimate.) The point is rather that *if* Gödel numbers are abused in a certain way, the proviso against ill-defined functions will be violated.

To illustrate: Suppose that $g(\varepsilon)$ is a function that takes in an expression ε and outputs its Gödel number. Then, define the function symbol '$r(x)$' as follows:

$$r(x) = \begin{cases} y & \text{if there is an } n \text{ such that } x = f_n(0) \text{ and } y = g(\ulcorner f_{\underline{n}}(0) \urcorner); \\ \uparrow & \text{otherwise.} \end{cases}$$

This defines '$r(x)$' in such a way that if, for some $n$, $x$ = the $0^{th}$ output of the $n$th function, then $r(x)$ = the Gödel number for the concatenation of the $n$th function symbol with '0'. So in particular, if '$x$' is replaced with an expression of the form $\ulcorner f_{\underline{n}}(0) \urcorner$, $r(x)$ then outputs the Gödel number for that self-same concatenation. In such a case, it outputs the Gödel number for the very sequence of symbols that replaces ‘$x$’. For example:

1. $r(f_0(0)) = g(`f_0(0)\text{'})$ [By definition of ‘$r(x)$’]
2. $r(f_1(0)) = g(`f_1(0)\text{'})$ [By definition of ‘$r(x)$’]

---

[10] Since any recursive function is numeralwise expressible in Peano Arithmetic, it may seem that that if contradiction is provable in the metalanguage, then contradiction is provable in the arithmetical object language. Yet the characteristic function for the set of reflection pairs $\langle a_n, \ulcorner a_{\underline{n}} \urcorner \rangle$ as defined below at (†) is not a well-defined function, much less a recursive one. So again, I doubt that formal arithmetic is in immediate danger.

And since each expression has its own Gödel number, we know that:

3. $g(`f_0(0)\text{'}) \neq g(`f_1(0)\text{'})$ [By definition of ‘$g(x)$’]

Suppose now that ‘$f_0(x)$’ and ‘$f_1(x)$’ have been earlier defined in such a way that:

4. $f_0(0) = f_1(0)$ [By definition of ‘$f_0(x)$’ and ‘$f_1(x)$’][11]

Then, given the truism that $f(f_0(0)) = f(f_0(0))$, the indiscernability of identicals assures us:

5. $r(f_0(0)) = r(f_1(0))$ [From 4]

However, the earlier observations entail:

6. $r(f_0(0)) \neq r(f_1(0))$ [From 1–3]

Yet lines 5 and 6 contradict.

Basically, ‘$r(x)$’ expresses a well-defined function only if no two function expressions of the form $\ulcorner f_{\underline{n}}(0) \urcorner$ co-refer. Hence, ‘$r(x)$’ does not express a well-defined function, for our language can be assumed to contain such co-referring expressions. The above definition of ‘$r(x)$’ must therefore be precluded. One might have thought, however, that it is acceptable to define a functional expression that denotes whichever Gödel number you please. But as our definition of ‘$r(x)$’ shows, if ‘$x$’ is replaced by any functional expression of the form $\ulcorner f_{\underline{n}}(0) \urcorner$, ‘$r(x)$’ *cannot* then be rigged to denote the Gödel number for that very functional expression. Otherwise, for some $\underline{m} \neq \underline{n}$, $\ulcorner f_{\underline{m}}(0) \urcorner$ will co-refer with $\ulcorner f_{\underline{n}}(0) \urcorner$, yet $\ulcorner r(f_{\underline{m}}(0)) \urcorner$ will denote a different Gödel number than $\ulcorner r(f_{\underline{n}}(0)) \urcorner$.

---

[11] E.g., 4 is true when ‘$f_0(x)$’ expresses the constantly 0 function and ‘$f_1(x)$’ expresses the identity function.

*5. Closing remark*

Thus, even Gödel numbering cannot be used without restriction. Yet as I say, since ill-defined functions are generally excluded, one can assume that existing formal systems are free of the problematic cases. Even so, it is worth acknowledging how such self-referential devices require *some* restraint, especially since the matter is glossed over by the usual "presumption of innocence" on self-reference.

*Appendix*

This appendix presents a different way of exploiting self-reference in a semantically open language to create an intensional functor. Notably, this avenue does not require quantification into subscript position; however, it instead requires an approach to defining function symbols that is not universally employed, although it is an approach which is commonly seen in Tarskian semantical treatments of formal languages.

The relevant definitional approach is one where a function is defined not directly but rather "indirectly" in virtue of defining a *symbol* for the function in the metalanguage. Thus, in a standard Tarskian approach, the function symbols are defined by "base clauses" (which are then used in connection with recursive clauses in defining the wffs of the object language). For example, if the domain is fixed as a domain of expressions, the Tarskian base clauses might include the following, where $\alpha$ is any singular term:

⌜$f(\alpha)$⌝ denotes ⌜'$\alpha$'⌝

Such a clause in the metalanguage introduces the quotation function for use in the object language: Given any expression, $f$ outputs the quotation of the expression; e.g., $f$('$a_0$') = ''$a_0$'' $f$('$a_1$') = ''$a_1$'', etc. (This is a well-defined function; it is not to be confused with the intensional

"quoting function" from section 1.) On some occasions, however, one might just simply introduce the function using '=' instead of 'denotes'. Where $x$ is any expression:

$f(x)$ = the quotation of $x$.

For the purposes of this appendix, however, the Tarskian "metalevel" approach is assumed.

The point, then, is that the possibilities for self-reference allow the introduction of the following intensional expression, where $\alpha$ is any singular term:

$\ulcorner i(\alpha)\urcorner$ denotes $\ulcorner i(\alpha)\urcorner$

For a given singular term $\alpha$, the concatenation of '*i*( )' with $\alpha$ is defined to denote that self-same concatenation. If self-reference were given free reign, such an expression would be thus far permissible. But suppose that a = b. Then, the definition implies:

*i*(a) = '*i*(a)'

*i*(b) = '*i*(b)'

Even though the input is the same, the output differs depending on which term for the input is used. Unrestrained self-reference is again sufficient by itself for an intensional function, and in this case, it is made possible without quantification into subscript position.